\documentclass[]{article}

 \usepackage[]{graphicx}
 \usepackage{epsfig}
 \usepackage{amsmath, amstext, amsfonts}
 \usepackage{amssymb}
 \usepackage{multicol}
 \usepackage[latin1]{inputenc}

\usepackage{anysize} % Soporte para el comando \marginsize
\marginsize{2cm}{2.5cm}{1.5cm}{2cm}

\begin{document}

 %.........................................................
 \title{Generalization of the Taylor formula}

 \author{D.V. Ionescu \\
 Department of Mathematics, University of Cluj, Cluj, Romania%\\ corradomassa@gmail.com
 %Aaron D Lanterman$^1$ and Jaemin Shin$^2$\footnote{To whom correspondence should be addressed.
  %InvPb-22-06.tex
  }

 \date{Gazeta Matematic\u{a} \c si Fizic\u{a} 7(60), pp. 389-395, August 1955 \\
%(online May 2009) pp. 30-32, Online ISBN: 9780511703768
%%\\{\small Received x, in final form x\\ Published 2005}
 }

 %%%% ===============================================================>>>>>>
 %%%% ===============================================================>>>>>>
 %%%% ===============================================================>>>>>>
 %%%% ===============================================================>>>>>>

%\begin{document}

\def\ct{\centerline}
\def\1{\'{\i}}
\def\={\!=\!}
\def\be{\begin{equation}}
\def\ee{\end{equation}}

\maketitle

%
%\begin{center}
%$^1$ Northumbria U, Scotland\\
%$^2$ Lancaster U., UK
%\end{center}

 %%%%%%%%%%%%%%%%%%%%%%%%%%%%%%% %%%%%%%%%%%%%%%%%%%%%%%%%%%
\begin{abstract}
\noindent The following expression
$$
y(x)=y_1(x)y(x_0)+y_2(x)y'(x_0)+\dots + y_n(x)y^{(n-1)}(x_0)+\int_{x_0}^{x}K(x,s)F[y(s)]ds
$$
is obtained as a generalization of the classical Taylor formula. The functions
$y_1(x),y_2(x),\dots, y_n(x)$ are the solutions of the differential equation $F(y) = y^{(n)}+a_1(x)y^{(n-1)}+\cdots +a_n(x)y=0$ fulfilling the initial conditions $y_{i}^{(k)}(x_0)=\delta _{ik}$, where $\delta_{ik}$ is the Kronecker symbol.
$K(x,s)$ is the solution of the equation $F(y)=0$ satisfying $K^{(i)}(x,s)_{_{x=s}}=\delta_{i,n-1}$.\\
{\tiny This paper has been originally published in Romanian. Translation by H.C. Rosu.}
%
%\noindent {\em Keywords}: Primordial black hole; Acoustic; Signal processing
\end{abstract}
%%%%%%%%%%%%%%%%%%%%%%%%%%%

\bigskip
\bigskip
\bigskip

%\section{Introduction}

\noindent {\bf I}. The theory of the adjoint differential equation to a given $n$th order linear differential equation \cite{dvi-1} allows to obtain a formula which generalizes the classical formula of the Taylor expansion:

\begin{equation}\label{e-1}
y(x) = y(x_0)+\frac{x-x_0}{1!}y'(x_0)+\cdots +\frac{(x-x_0)^{n-1}}{(n-1)!}y^{(n-1)}(x_0)+\int _{x_0}^{x}\frac{(x-s)^{n-1}}{(n-1)!}y^{(n)}(s)ds~.
\end{equation}

Let
\begin{equation}\label{e-2}
F(y) = y^{(n)}+a_1(x)y^{(n-1)}+\cdots +a_n(x)y=0
\end{equation}
be a homogeneous linear differential equation of the $n$th order whose coefficients $a_i(x)$, $i=1,2,...,n$, are continuous functions on the interval $[a,b]$.

Let also

\begin{equation}\label{e-3}
G(z)=(-1)^n[z^{(n)}-(a_1z)^{(n-1)}+\cdots +(-1)^na_nz]=0
\end{equation}
be the adjoint equation of eq.~(\ref{e-2}).

\medskip

The following identity is known to be true

\begin{equation}\label{e-4}
zF(y)-yG(z) = \frac{d}{dx}\Upsilon(y,z)~,
\end{equation}
where
\begin{align}\label{e-5}
\Upsilon(y,z)=&zy^{(n-1)}+[(a_1z)-z']y^{(n-2)}+[(a_2z)-(a_1z)'+z'']y^{(n-3)}+\cdots \nonumber\\
&+[(a_{n-1}z)-(a_{n-2}z)'+\cdots +(-1)^{n-1}z^{(n-1)}]y~.
\end{align}

Let us integrate both sides of (\ref{e-4}) from $x_0$ to $x$, where both of them are points belonging to the interval $[a,b]$. We get
\begin{equation}\label{e-6}
\int _{x_0}^{x}\{ z(s)F(y(s))-y(s)G(z(s))\}ds =\int _{x_0}^{x} \frac{d}{ds}\Upsilon(y(s),z(s))ds~.
\end{equation}
Let
\begin{equation}\label{e-7}
z(s)=\varphi(x,s)
\end{equation}
be a solution of the adjoint equation $G(z)=0$ which depends on the parameter $x$ and fulfills the initial conditions
\begin{align}\label{e-8}
&z(s)|_x=\varphi(x,x)=0~, \, z'(s)|_x=\varphi'_{s}(x,x)=0~,\dots~, \, z^{(n-2)}(s)|_x=\varphi^{(n-2)}_{s}(x,x)=0~, \nonumber \\
&z^{(n-1)}(s)|_x=\varphi^{(n-1)}_{s}(x,x)=1~.
\end{align}
We now replace $z(s)$ by $\varphi(x,s)$ in eq.~(\ref{e-6}) which leads to the following equation:
\begin{equation}\label{e-9}
\int_{x_0}^{x}\varphi(x,s)F[y(s)]ds=\Upsilon[y(s),z(s)]|_{x_0}^{x}~.
\end{equation}
It is worth noticing that $\Upsilon(y,z)$ is a linear and homogeneous function in $z,z',\dots,z^{(n-1)}$, the coefficient of $z^{(n-1)}$ being
$(-1)^{n-1}y$. For $x=s$, $\Upsilon(y,z)$ reduces to $(-1)^{n-1}y(s)$ because of (\ref{e-8}). On the other hand, $\Upsilon(y,z)$ is a linear and homogeneous function in $y,y',\dots,y^{(n-1)}$ whose coefficients are displayed in (\ref{e-5}). They depend on $x$ and $s$ because $z=\varphi(x,s)$.
For $s=x_0$, the coefficients of $y(x_0),y'(x_0),\dots,y^{(n-1)}(x_0)$ are functions of $x$ that we denote by
$$
(-1)^{n-1}\lambda_1(x)~,\, (-1)^{n-1}\lambda_2(x)~,\,\dots,\, (-1)^{n-1}\lambda_n(x)~.
$$
Equation (\ref{e-9}) becomes
$$
\int_{x_0}^{x}\varphi(x,s)F[y(s)]ds=(-1)^{n-1}y(x)+(-1)^{n-1}[\lambda_1(x)y(x_0)+\lambda_2(x)y'(x_0)+\dots + \lambda_n(x)y^{(n-1)}(x_0)]~.
$$
Solving for $y(x)$, one gets
\begin{equation}\label{e-10}
y(x)=\lambda_1(x)y(x_0)+\lambda_2(x)y'(x_0)+\dots +\lambda_n(x)y^{(n-1)}(x_0)+(-1)^{n-1}\int_{x_0}^{x}\varphi(x,s)F[y(s)]ds~.
\end{equation}
Let us give the meaning of the functions $\lambda_i(x)$. For this, in eq.~(\ref{e-10}) we substitute $y(x)$ by the solution $y_i(x)$ of the differential equation $F(y)=0$ which satisfies the initial conditions
\begin{equation}\label{e-11}
y_i(x_0)=0,\dots ,\,y_{i}^{(i-2)}(x_0)=0,\, y_{i}^{(i-1)}(x_0)=1,\, y_{i}^{(i)}(x_0)=0,\dots ,\,y_{i}^{(n-1)}(x_0)=0,\,
\end{equation}
where $i=1,2,\dots,n$. Thus, we obtain
$$
y_i(x)=\lambda_i(x)~.
$$
We conclude that $\lambda_1(x),\lambda_2(x),\dots,\lambda_n(x)$ in eq.~(\ref{e-10}) are the solutions of the differential equation $F(y)=0$ with the initial conditions shown in the following Table:

\medskip

\begin{center}
$%
\begin{tabular}
[b]{|c|c|}\hline
%\multicolumn{2}{|c|}{Mathieu's $\lambda_{n}$}%
%\\\hline
$y_1$ & $1,0,\dots,0$ \\\hline
$y_2$ & $0,1,\dots,0$ \\\hline
$\cdots$ & $\cdots\cdots\cdots$ \\\hline
$y_n$ & $0,0, \dots,1$ \\\hline
\end{tabular}
\qquad$%
\end{center}

\medskip

It is known that in this case $y_1(x),y_2(x),\dots, y_n(x)$ is a fundamental set of solutions of the differential equation $F(y)=0$.
Thus, equation (\ref{e-10}) turns into
\begin{equation}\label{e-13}
y(x)=y_1(x)y(x_0)+y_2(x)y'(x_0)+\dots +y_n(x)y^{(n-1)}(x_0)+(-1)^{(n-1)}\int_{x_0}^{x}\varphi(x,s)F[y(s)]ds~.
\end{equation}

\medskip

In the particular case when
$$
F(y)=y^{(n)}~,
$$
one obtains
$$
G(z)=(-1)^{n}z^{(n)}~.
$$
The solution of the differential equation $G(z)=0$ satisfying the initial conditions (\ref{e-8}) is
$$
\varphi(x,s)=(-1)^{n-1}\frac{(x-s)^{n-1}}{(n-1)!}
$$
and the solutions of the differential equation $y^{(n)}=0$ which satisfy the initial conditions given for each of them in the Table are
\begin{align}%\label{e-14}
y_1(x)&=1\nonumber\\
y_2(x)&=\frac{x-x_0}{1!}\nonumber\\
\dots & \dots\dots\dots\dots\nonumber\\
y_n(x)&=\frac{(x-x_0)^{n-1}}{(n-1)!}~.\nonumber
\end{align}
Formula (\ref{e-13}) turns in this case into Taylor's expansion in (\ref{e-1}). This explains why the formula (\ref{e-13}) on which we want to draw the attention of the reader of this work is a generalization of the Taylor formula.

\bigskip

\noindent {\bf II}. The method just described allows to express the fundamental system of solutions $y_1(x),y_2(x),\dots, y_n(x)$ of the differential equation $F(y)$ which satisfy the initial conditions written for each of them in the Table above by means of the functions $\phi(x,s)$.

Indeed, taking into account the expression (\ref{e-4}) of $\Upsilon(y,z)$, we find out that
 \begin{align}\label{e-14}
 y_1(x)&=(-1)^{n}\bigg[a_{n-1}(s)\phi(x,s)-\frac{\partial(a_{n-2}(s)\phi(x,s))}{\partial s}+\dots +(-1)^{n-1}\frac{\partial^{n-1}\varphi(x,s)}{\partial s^{n-1}} \bigg]_{s=x_0}\nonumber \\
 \dots &\dots \dots \dots \dots \dots \dots \dots\dots\dots\dots\dots\dots\dots\dots\dots\dots\dots\dots\dots\dots\\
 y_{n-2}(x)&=(-1)^{n}\bigg[a_{2}(s)\phi(x,s)-\frac{\partial(a_{1}(s)\phi(x,s))}{\partial s} +\frac{\partial^{2}\varphi(x,s)}{\partial s^{2}} \bigg]_{s=x_0}\nonumber \\
y_{n-1}(x)&=(-1)^{n}\bigg[a_{1}(s)\phi(x,s) -\frac{\partial\varphi(x,s)}{\partial s} \bigg]_{s=x_0}\nonumber \\
 y_n(x)& =(-1)^{(n-1)}\varphi(x,s)|_{s=x_0} \nonumber
\end{align}

\medskip

Focusing our attention on the last formula of this set, we infer that if $s$ is any point in $[a,b]$, the integral $K(x,s)$ of $F(y)=0$ which depends on the parameter $s$, and for which the initial conditions
\begin{equation}\label{e-15}
K(x,s)|_{_{x=s}}=0~, \quad K'(x,s)|_{_{x=s}}=0~,\dots K^{(n-2)}(x,s)|_{_{x=s}}=0~, \quad K^{(n-1)}(x,s)|_{_{x=s}}=1~,
\end{equation}
hold, is $(-1)^{(n-1)}\varphi(x,s)$. Thus, we have the identity
\begin{equation}\label{e-16}
K(x,s)=(-1)^{(n-1)}\varphi(x,s)~.
\end{equation}

\medskip

One consequence of this identity is that the fundamental formula (\ref{e-13}) can also be written as follows
\begin{equation}\label{e-13bis}
y(x)=y_1(x)y(x_0)+y_2(x)y'(x_0)+\dots + y_n(x)y^{(n-1)}(x_0)+\int_{x_0}^{x}K(x,s)F[y(s)]ds~.
\end{equation}

\bigskip

\noindent {\bf III}. Applications\\

1. Let us consider the nonhomogeneous differential equation
$$
F(y)=f(x)~.
$$
Its solution for which the following initial conditions
$$
y(x_0)=0,\,\, y'(x_0)=0,\dots , y^{(n-1)}(0)=0
$$
hold is given by eq.~(\ref{e-13bis}) where in the right-hand side $F[y(s)]$ is substituted by $f(s)$. One gets in this way

\begin{equation}\label{e-17}
y(x) = \int _{x_0}^{x}K(x,s)f(s)ds~,
\end{equation}
which is the well-known Cauchy formula.

\bigskip

2. Assume that
$$
F(y)=y''+y~.
$$
We will have
$$
G(z)=z''+z~
$$
and the solution of the equation
$$
z''+z=0
$$
which fulfills the conditions
$$
z(s)|_{_{s=x}}=0~, \quad z's)|_{_{s=x}}=1~,
$$
is
$$
z=\varphi(x,s)=-\sin(x-s)~.
$$
Equations (\ref{e-14}) provide
$$
y_1(x)=\cos(x-x_0)~, \qquad y_2(x)=\sin(x-x_0)
$$
and eq.~(\ref{e-13}) becomes
$$
y(x)=y(x_0)\cos(x-x_0)+y'(x_0)\sin(x-x_0)+\int_{x_0}^{x}\sin(x-s)[y''(s)+y(s)]ds~.
$$
The following formula:
$$
y(x)=y(x_0)\cosh(x-x_0)+y'(x_0)\sinh(x-x_0)+\int_{x_0}^{x}\sinh(x-s)[y''(s)-y(s)]ds
$$
can be similarly proved.

\medskip

3. Let us take
$$
F(y)=y''''+5y''+4y~,
$$
which implies:
$$
G(z)=z''''+5z''+4z~.
$$
The solution of $G(z)=0$ which fulfills the initial conditions
$$
z(s)|_{_{s=x}}=0~, \quad z'(s)|_{_{s=x}}=0~,\quad z''(s)|_{_{s=x}}=0~, \quad z'''(s)|_{_{s=x}}=1
$$
is
$$
z=\varphi(x,s)=\frac{\sin 2(x-s)-2\sin(x-s)}{6}~.
$$
Applying equations (\ref{e-14}), one finds
\begin{align*}
y_1(x)&=\frac{1}{3}\bigg[4\cos(x-x_0)-\cos2(x-x_0)\bigg]~, \nonumber\\
y_2(x)&=\frac{1}{6}\bigg[8\sin(x-x_0)-\sin 2(x-x_0)\bigg]~, \nonumber\\
y_3(x)&=\frac{1}{3}\bigg[\cos(x-x_0)-\cos2(x-x_0)\bigg]~, \nonumber\\
y_4(x)&=\frac{1}{6}\big[2\sin(x-x_0)-\sin 2(x-x_0)\bigg]~. \nonumber\\
\end{align*}
Formula (\ref{e-13}) becomes
\begin{align*}\label{e-19}
y(x)&=\frac{1}{3}\big[4\cos(x-x_0)-\cos2(x-x_0)\big]y(x_0)+\frac{1}{6}\big[8\sin(x-x_0)-\sin 2(x-x_0)\big]y'(x_0)\nonumber\\
&+ \frac{1}{3}\big[\cos(x-x_0)-\cos2(x-x_0)\big]y''(x_0)+\frac{1}{6}\big[2\sin(x-x_0)-\sin 2(x-x_0)\big]y'''(x_0)\nonumber\\
&-\int _{x_0}^{x}\frac{1}{6}\big[\sin 2(x-s)-2\sin(x-s)\big]\big[y''''(s)+5y''(s)+4y(s)\big]ds~.\nonumber\\
\end{align*}

\medskip

4. Let us consider the integro-differential equation
\begin{equation}\label{e-20}
y^{(n)}+a_1(x)y^{(n-1)}+\ldots +a_n(x)y=\int _{x_0}^{x}N(x,s)y(s)ds+f(x)~,
\end{equation}
where $a_1(x),\dots, a_n(x)$ are continuous functions in the interval $[a,b]$, $x_0\in [a,b]$, and $N(x,s)$ is a continuous function in the domain
$a\leq x\leq b$, $a\leq s\leq b$. We want to show that the solution of (\ref{e-20}) fulfilling the initial conditions
\begin{equation}\label{e-21}
y(x_0)=y_0~, \quad y'(x_0)=y'_0~, \quad \dots, \quad y^{(n-1)}(x_0)=y_{0}^{(n-1)}
\end{equation}
is the solution of an integral Volterra equation of the second kind which we will obtain in the following.
Let us apply eq.~(\ref{e-13bis}), where $F(y)$ is the right hand side of eq.~(\ref{e-20}). Then, we get
\begin{equation}\label{e-22}
y(x)=y_1(x)y_0+y_2(x)y'_0+\dots +y_n(x)y_{0}^{(n-1)}+ \int _{x_0}^{x}K(x,s)ds\int _{x_0}^{s}N(x,t)y(t)dt+\int _{x_0}^{x}K(x,s)f(s)ds~.
\end{equation}
Let $Y(x)$ be the solution of the differential equation $F(y)=f(x)$ which satisfies the initial conditions (\ref{e-21}). Its expression is given by
\begin{equation}\label{e-23}
Y(x)=y_1(x)y_0+y_2(x)y'_0+\dots +y_n(x)y_{0}^{(n-1)}+ \int _{x_0}^{x}K(x,s)f(s)ds~.
\end{equation}
On the other hand, we notice that
\begin{equation}\label{e-24}
N_1(x,t)=\int _{t}^{x}K(x,s)N(s,t)dt~.
\end{equation}
Therefore, eq.~(\ref{e-22}) becomes
\begin{equation}\label{e-25}
y(x)=\int _{x_0}^{x}N_1(x,t)y(t)dt+Y(x)~.
\end{equation}
Thus, we have shown that the solution of the integro-differential equation (\ref{e-20}) which fulfills the initial conditions (\ref{e-21}) is the solution of the integral Volterra equation of the second kind with the kernel $N_1(x,t)$ given by (\ref{e-24}), while $Y(x)$ is the solution of the differential equation $F(y)=f(x)$ with the same initial conditions (\ref{e-21}).

%eqs.~(\ref{e-3}) and (\ref{e-4})

\bigskip
\bigskip
\bigskip

%
%We would like to thank Dr. X for a careful reading of the first draft of this work.
%The second author wishes to thank Y for partial support through project 46980.
%
%

\end{document}